\documentclass[]{amsart}
\usepackage{graphicx}
\usepackage{amsmath}

\newtheorem{defn}{Definition}
\newtheorem{rem}{Remark}
\newtheorem{exa}{Example}

\begin{document}

\title[Topologies of Non-intersecting Circles in the Plane]{Topologically Distinct Sets of Non-intersecting Circles in the Plane}

\author{Richard J. Mathar}
\urladdr{http://www.mpia.de/~mathar}
\address{Max-Planck Institute of Astronomy, K\"onigstuhl 17, 69117 Heidelberg, Germany}

\date{\today}

\begin{abstract}
Nested parentheses are forms in an algebra which define orders of evaluations. A class of well-formed
sets of associated opening and closing parentheses is well studied in conjunction
with Dyck paths and Catalan numbers. Nested parentheses also represent cuts through
circles on a line. These become topologies of non-intersecting circles in the plane
if the underlying algebra is commutative. 

This paper generalizes the concept and answers quantitatively---as recurrences
and generating functions of matching rooted forests---the questions: how many
different topologies of nested circles exist in the plane if (i) pairs of circles
may intersect, or (ii) even triples of circles may intersect. That analysis is driven
by examining the symmetry properties of the inner regions of the fundamental type(s)
of the intersecting pairs and triples.
\end{abstract}

\maketitle

\section{Paired Parentheses and Catalan Numbers}\label{sec.P}
In a (non-commutative) algebra, opening and closing
parentheses prescribe the order of grouping and
evaluating expressions. 
\begin{defn}
A string of parentheses is \emph{well-formed}
if the total number of opening parentheses equals the number
of closing parentheses, and if the subtotal count of opening parentheses
is always larger than or equal to the subtotal count of closing parentheses 
while parsing the string left-to-right. 
\end{defn}
\begin{rem}
Equivalently one may demand that the subtotal of closing parentheses
is always larger or equal to the subtotal of opening parentheses
while parsing the string right-to-left.
\end{rem}
The well-formed nested parentheses form sets $\mathbb{P}_N$
of expressions with $N$ pairs of parentheses. 

\begin{defn}
$\mathbb{P}_N$ is the set of all well-formed expressions with $N$
pairs of parentheses.
\end{defn}
There are expressions
that can be \emph{factored}---in the algebra---by cutting the string
at some places such that the left and right substrings are also well-formed.
The number $f$ of their factors puts the elements of $\mathbb{P}_N$
into disjoint subsets:

\begin{defn}
$\mathbb{P}_N^{(f)}$ is the set of all well-formed expressions with $N$
pairs of parentheses and $1\le f\le N$ factors.
\end{defn}
\begin{eqnarray}
\mathbb{P}_N &=& \bigcup_f \mathbb{P}_N^{(f)};\\
|\mathbb{P}_N| &=& \sum_{f=1}^N |\mathbb{P}_N^{(f)}|.
\end{eqnarray}
\begin{exa}
\begin{eqnarray}
\mathbb{P}_1 &=& \mathbb{P}_1^{(1)}=\{()\}; \\
\mathbb{P}_2^{(1)} &=& \{(())\}; \\
\mathbb{P}_2^{(2)} &=& \{()()\}; \\
\mathbb{P}_3^{(1)} &=& \{((())), (()())\}; \\
\mathbb{P}_3^{(2)} &=& \{(())(), ()(())\}; \\
\mathbb{P}_3^{(3)} &=& \{()()()\}; \\
\mathbb{P}_4^{(1)} &=& \{(((()))), ((()())), ((())()), (()(())), (()()())\};
\end{eqnarray}
\end{exa}

\begin{rem}
The opening and closing parentheses are the two letters in an alphabet of
words, with a grammar that recursively admits words 
\begin{enumerate}
\item
that are the empty word,
\item
that are concatenations of two words,
\item
that are concatenations of the first letter, a word, and the second letter.
\end{enumerate}
\end{rem}

If the opening parentheses are replaced by U and the closing parentheses
replaced by D an equivalence with Dyck paths arises; the number of
returns to the horizontal line in the paths is equivalent to the
number of factors in the expression. This leads straight to the
well-known Catalan triangle \cite[A033184]{EIS} of Table \ref{tab.P}.
The set sizes $|\mathbb{P}_N|$ are 
the Catalan numbers 
\cite[\S 1.15]{Comtet}\cite[A000108]{EIS}\cite{DonagheyJCTB22}.

\begin{table}
\begin{tabular}{r|r|rrrrrrrrrr}
$N$ & $ |\mathbb{P}_N|$ & 1 & 2 & 3 & 4 & 5 & 6 & 7 &8 &9\\
\hline
1 & 1 & 1 \\
2 & 2 & 1 & 1\\
3 & 5 & 2 & 2 & 1 \\
4 & 14 & 5 & 5 & 3 & 1\\
5 & 42 & 14 & 14 & 9 & 4 & 1\\
6 & 132 &42&42&28&14&5&1 \\
7 & 429 &132&132&90&48&20&6&1 \\
8 & 1430 &429&429&297&165&75&27&7&1 \\
9 & 4862 &1430&1430&1001&572&275&110&35&8&1 \\
\end{tabular}
\caption{Catalan triangle: The number of nested expressions with $N$ pairs
of parentheses: the total count $|\mathbb{P}_N|$ and the
number of nested expressions with $1\le f\le N$ factors, $|\mathbb{P}^{(f)}_N|$.
}
\label{tab.P}
\end{table}

\begin{rem}
Because an expression distributes $N$ left parenthesis at $2N$ places,
the set size is limited by $|\mathbb{P}_N|<\binom{2N}{N}$, the central
binomial coefficients. 
Actually one of them must be placed at the leftmost place and 
none can be placed at the rightmost place, which
leaves $2N-2$ places to distribute $N-1$ of them:
$|\mathbb{P}_N|<\binom{2N-2}{N-1}$.
\end{rem}

\begin{rem}\label{rem.bits}
A computer representation uses the two binary digits
\texttt{1} and \texttt{0} to represent the opening and the closing
parenthesis in the aforementioned alphabet of two letters
\cite[A063171,A014486]{EIS}.
(Then the most-significant bit is always 1\@. 
The swapped mapping is less useful because it needs to deal
with the numerical representation of leading zeros in the binary number.)
Because the rightmost part is absent for the unique
case of missing parentheses, $N=0$, or a closing parenthesis, all
these representations are even numbers.
This representation by non-negative integers induces a strict
ordering in the set of nested parentheses.
\begin{exa}
$() = 10\_2$;
$()() = 1010\_2$;
$(()) = 1100\_2$;
$()()() = 101010\_2$ ;
$()(()) = 101100\_2$ ;
$(())() = 110010\_2$ ;
$(()()) = 110100\_2$ ;
$((())) = 111000\_2$.
\end{exa}
\end{rem}

The expressions with one factor are given by
embracing any expression with one pair less
at the left and right end with a pair of matching parentheses:
\begin{equation}
|\mathbb{P}_N^{(1)}| = |\mathbb{P}_{N-1}|.
\label{eq.Pf1}
\end{equation}

The number of expressions with $f$ factors
is given by considering any concatenated ``word'' of 
factorizations \cite{KlarnerJCT9},
\begin{equation}
|\mathbb{P}_N^{(f)}| = \sum_{C(N): N=N_1+N_2+\cdots+N_f} 
|\mathbb{P}_{N_1}^{(1)}|
|\mathbb{P}_{N_2}^{(1)}|
\cdots
|\mathbb{P}_{N_f}^{(1)}|,
\quad f\ge 2,
\label{eq.Pf}
\end{equation}
where the sum is over all compositions (``ordered'' partitions)
of $N$ into positive parts $N_j$ such that subexpressions 
do not factor any further.

A well-formed expression of parentheses represents a 
set of $N$ nested circles if we join the upper and lower end of
each associated pair of parentheses. The radii of the circles
are growing functions of their spatial distance in the expressions; their
mid points are on a
straight line, and no perimeters of
any pair of circles intersect. The string of opening and closing
parentheses is a record of entering or leaving a circle
while poking from the outside along the line through all circles.
For each pair of circles (i) either the smaller one is entirely
immersed in the larger one or (ii) they have no common points.
\setlength{\unitlength}{1cm}
\begin{exa}\label{exa.2d}
$(()())() \mapsto$
\begin{picture}(7,1.5)
\put(1.0,0.5){\circle{0.5}}
\put(1.7,0.5){\circle{0.5}}
\put(1.35,0.5){\circle{3.0}}
\put(2.5,0.5){\circle{0.5}}
\end{picture}
\end{exa}

\section{Nonintersecting Circle Sets in the Plane}\label{sec.C} 
\subsection{Nested Circle Sets}
If the algebra of Section \ref{sec.P} is
a commutative algebra, some sets of nested parentheses
are no longer considered distinct, because the order
of the factors does no longer matter.
\begin{defn}
$\mathbb{C}_N$ is the set of well-formed expressions
of $N$ pairs of parentheses where the order within factorizations
does not matter.
\end{defn}
\begin{defn}
$\mathbb{C}_N^{(f)}$ is the set of well-formed expressions
of $N$ pairs of parentheses with $f$ factors
where the order within factorizations
does not matter.
\end{defn}
The number of factors still is a unique parameter of each well-formed
set of expressions, so
\begin{eqnarray}
\mathbb{C}_N &=& \bigcup_f \mathbb{C}_N^{(f)};\\
|\mathbb{C}_N| &=& \sum_{f=1}^N |\mathbb{C}_N^{(f)}|;\quad |\mathbb{C}_0|=1.
\label{eq.Cunion}
\end{eqnarray}
We ``lose'' some of the sets of parenthesis relative to Section \ref{sec.P},
because for example now the expressions $()(())$ and $(())()$ are considered
the same: $|\mathbb{C}_N|\le |\mathbb{P}_N|$. This reduction in the
admitted expressions applies recursively to all sub-expressions that
are obtained by ``peeling'' the surrounding pair of parentheses off
expressions with a single factor.
\begin{rem}
The reduction of equivalent expressions to a single representation
requires some convention of which ordering of the factors is the
admitted one. One convention is to map each factor to an integer
with the representation of Remark \ref{rem.bits}, to put these factors
into non-increasing or non-decreasing numerical order, and to 
concatenate their binary representations left-to-right to define
the representative.
\end{rem}
\begin{exa}
\begin{eqnarray}
\mathbb{C}_1 &=& \mathbb{C}_1^{(1)}=\{()\}; \\
\mathbb{C}_2^{(1)} &=& \{(())\}; \\
\mathbb{C}_2^{(2)} &=& \{()()\}; \\
\mathbb{C}_3^{(1)} &=& \{((())), (()())\}; \\
\mathbb{C}_3^{(2)} &=& \{(())()\}; \\
\mathbb{C}_3^{(3)} &=& \{()()()\}; \\
\mathbb{C}_4^{(1)} &=& \{(((()))), ((()())), ((())()), (()()())\};
\end{eqnarray}

\end{exa}

Table \ref{tab.C} shows $\mathbb{C}_N^{(f)}$ and their sums $|\mathbb{C}_N|$.
The values are bootstrapped as follows:
\begin{table}
\begin{tabular}{r|r|rrrrrrrrrrrrrrrrrrrrr}
$N$ & $ |\mathbb{C}_N|$ & 1 & 2 & 3 & 4 & 5 & 6 & 7 &8 & 9 & 10 & 11\\
\hline
        1 & 1&1 \\
        2 & 2&1&1 \\
        3 & 4&2&1&1 \\
        4 & 9&4&3&1&1 \\
        5 & 20&9&6&3&1&1 \\
        6 & 48&20&16&7&3&1&1 \\
        7 & 115&48&37&18&7&3&1&1 \\
        8 & 286&115&96&44&19&7&3&1&1 \\
        9 & 719&286&239&117&46&19&7&3&1&1 \\
        10 & 1842&719&622&299&124&47&19&7&3&1&1 \\
        11 & 4766&1842&1607&793&320&126&47&19&7&3&1&1 \\
        12 & 12486&4766&4235&2095&858&327&127&47&19&7&3&1&1 \\
\end{tabular}
\caption{The counts $|\mathbb{C}_N|$ and
$|\mathbb{C}_N^{(f)}|$ of nested nonintersecting circles \cite[A033185,A000081]{EIS}.
}
\label{tab.C}
\end{table}
The consideration leading to Equation (\ref{eq.Pf1}) leads also to
\begin{equation}
|\mathbb{C}_N^{(1)}| =
|\mathbb{C}_{N-1}|.
\label{eq.CN1}
\end{equation}
The decomposition of an expression with $f$ factors needs
to consider the number of ways of distributing the $N$ pairs of
parentheses over elements that do not factorize further. We partition
$N$ as $\pi(N): N=\{N_1^{c_1};N_2^{c_2};\ldots N_f^{c_f}\}=\sum_{j=1}^f {c_jN_j}$ meaning
that the expression contains $c_1$ factors with elements of 
$\mathbb{C}_1^{(1)}$,
$c_2$ factors with elements of $\mathbb{C}_2^{(1)}$, and so on.
For each part $N_j$ with repetition $c_j$ we compute 
the number of lists  of $c_j$
elements taken from a set of $|\mathbb{C}_{N_j}^{(1)}|$, possibly
selecting some elements more than once or not at all.
This is the number of weak compositions of $c_j$ into $C_{N_j}^{(1)}$
parts of non-negative integers.
This equals the number of compositions
of $c_j+|\mathbb{C}_{N_j}^{(1)}|$ into $|\mathbb{C}_{N_j}^{(1)}|$ parts of positive integers,
which is $\binom{c_j+|\mathbb{C}_{N_j}^{(1)}|-1}{c_j}$ \cite[\S 1.2]{Stanley}.
\begin{equation}
|\mathbb{C}_N^{(f)}| = \sum_{\pi(N): N=\{N_1^{c_1};N_2^{c_2};\ldots N_f^{c_f}\}} 
\binom{|\mathbb{C}_{N_1}^{(1)}|+c_1-1}{c_1}
\binom{|\mathbb{C}_{N_2}^{(1)}|+c_2-1}{c_2}
\cdots
\binom{|\mathbb{C}_{N_f}^{(1)}|+c_f-1}{c_f}.
\label{eq.CNf}
\end{equation}

Table \ref{tab.C} shows the phenomenon that at large $N$
the values at large $f$ converge to the sequence
\begin{equation}
|\mathbb{C}_N^{(N-i)}|=1,1,3,7,19,47,127,330,889,2378,\ldots, \quad i\ge 0,\quad N\to \infty,
\end{equation}
the \emph{envelope} as Knopfmacher and Mays call it \cite{KnopfmacherAC53}.
This is the Euler transform of $|\mathbb{C}_N|$ \cite{AgnewAJM66} and means
that if the number of factors is large, most of the factors are the element $\mathbb{C}_1{(1)}=\{()\}$
and only few combinations remain to exhaust the others.

Returning to the interpretation of $\mathbb{P}_N$ as
non-intersecting circles on a line, considering the order
of factorizations unimportant means that $\mathbb{C}_N$
contains topologically distinct sets of non-intersecting
circles that are free to move away from the line---as long as they
stay within the boundaries of their surrounding circles.
The two circles inside the bigger circle in Example \ref{exa.2d} are
allowed to bump around within the bigger circle, and the big and 
the outer small circle may also jointly move to other places.
\begin{rem}
This is a planetary model of the circles in the sense that each circle
can ``rotate'' inside its surrounding circle, and all
these geometries are considered equivalent.
\end{rem}

\begin{defn} The generating function for the number of nested
expressions in the commutative algebra is
\begin{equation}
C(z)=\sum_{N\ge 0} |\mathbb{C}_N|z^N.
\label{eq.Cz}
\end{equation}
\end{defn}
It satisfies \cite[I.5.2]{Flajolet}\cite{PolyaAA68,CayleyPMag4}
\begin{equation}
C(z)= \exp\left(\sum_{j\ge 1} z^jC(z^j)/j\right).
\end{equation}

\subsection{Nested Circles Embedded in the Sphere Surface}
If the $N$ circles are not embedded in the plane but embedded in the
surface of a three-dimensional sphere, the topologies are counted
by the unlabeled trees with $N+1$ nodes as stated
by Reshetnikov \cite[A000055]{EIS}:
\begin{equation}
1, 1, 1, 2, 3, 6, 11,23,47,106,\ldots\quad N\ge 0.
\end{equation}
Each unlabeled tree can be mapped to a circle set topology
by constructing the line graph of the tree, associating
circles with nodes of the line graph,  and arranging
the circles on the sphere such that they can touch (by moving them on the sphere
and changing radii)
iff they are connected in the line graph. The tree is
a connectivity diagram of the regions on the sphere surface: an edge
in the tree indicates one must cross a circle boundary to enter a different region.

\begin{figure}
\includegraphics[scale=0.5]{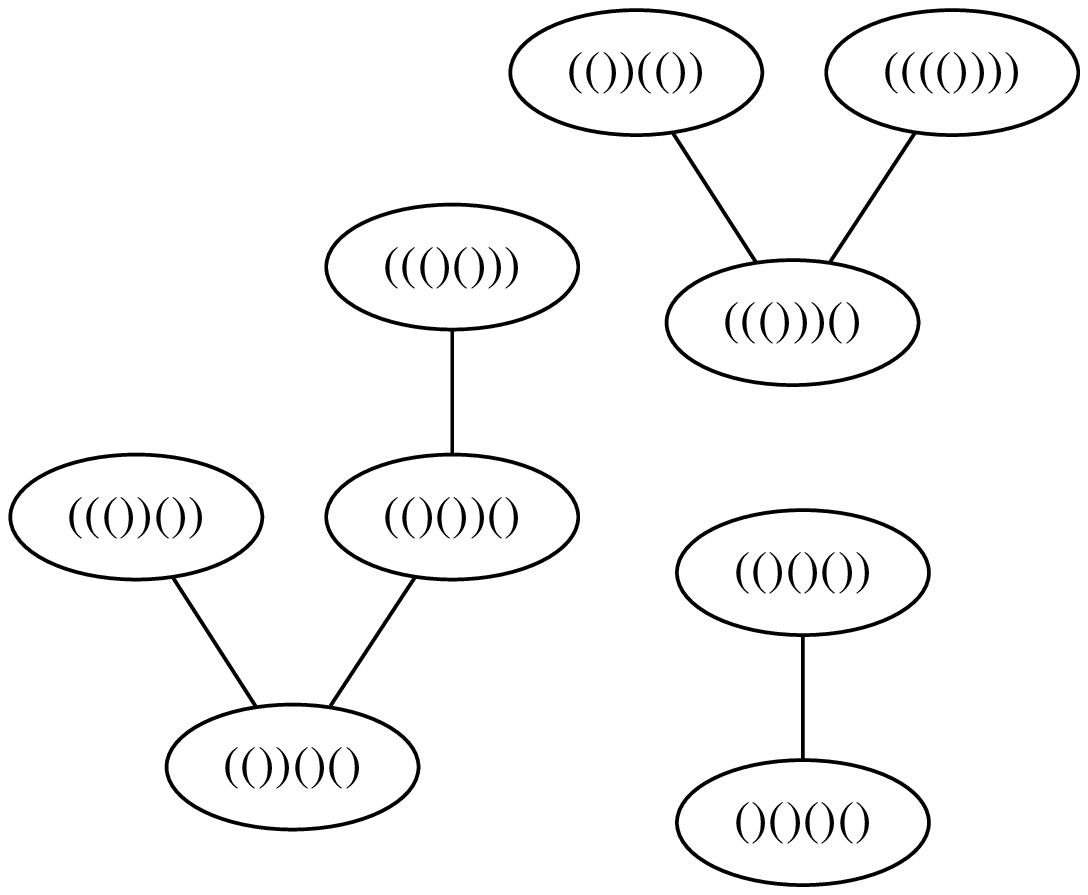}
\caption{The 3 clusters of grouping the $|\mathbb{C}_4|=9$
expressions with 4 pairs of parentheses into clusters of expressions equivalent under the flip transform.
}
\label{fig.C4}
\end{figure}

\begin{figure}
\includegraphics[scale=0.5]{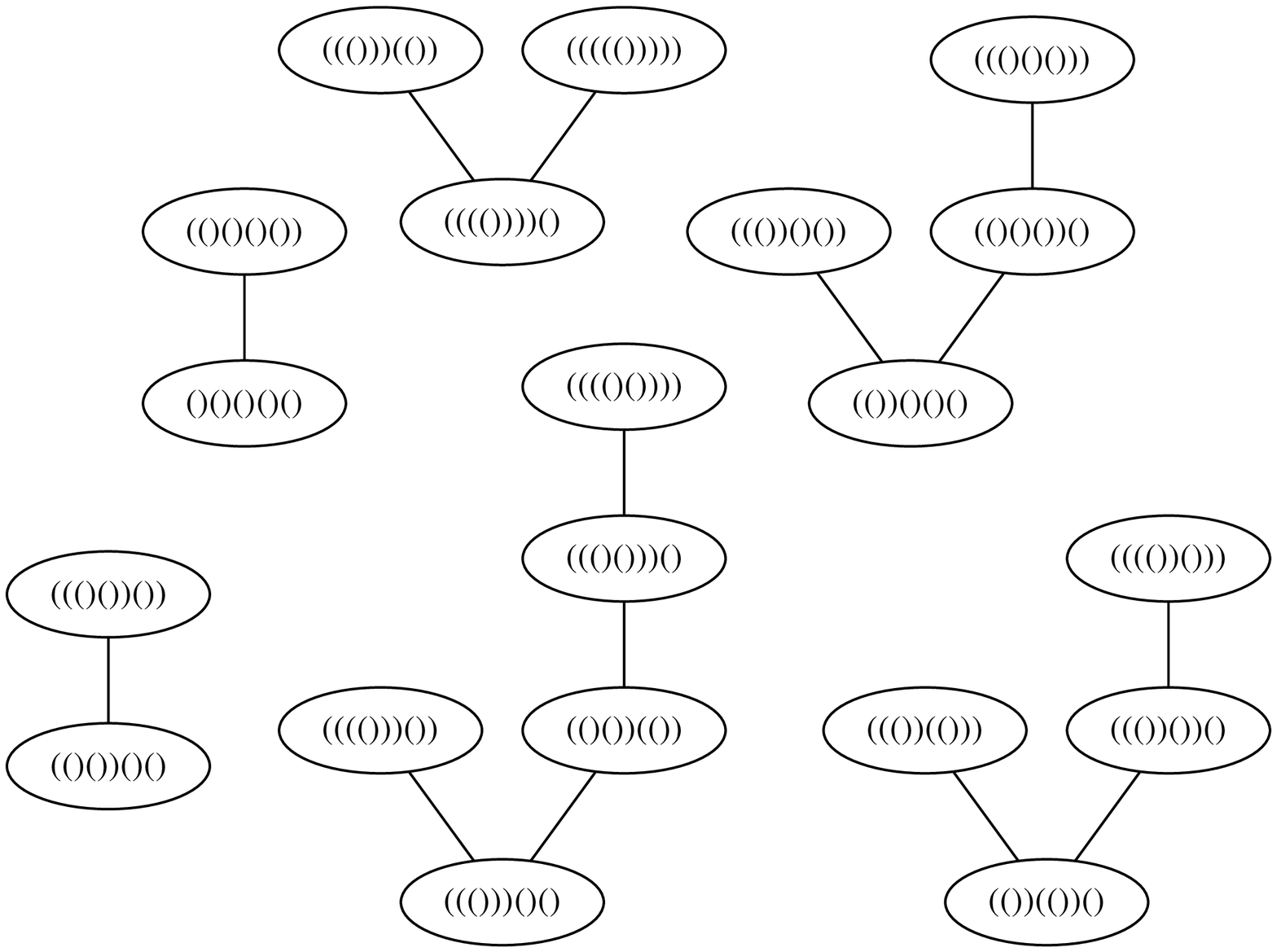}
\caption{The 6 clusters of grouping the $|\mathbb{C}_5|=20$
expressions with 5 pairs of parentheses into clusters.
}
\label{fig.C5}
\end{figure}

\begin{figure}
\includegraphics[scale=0.5]{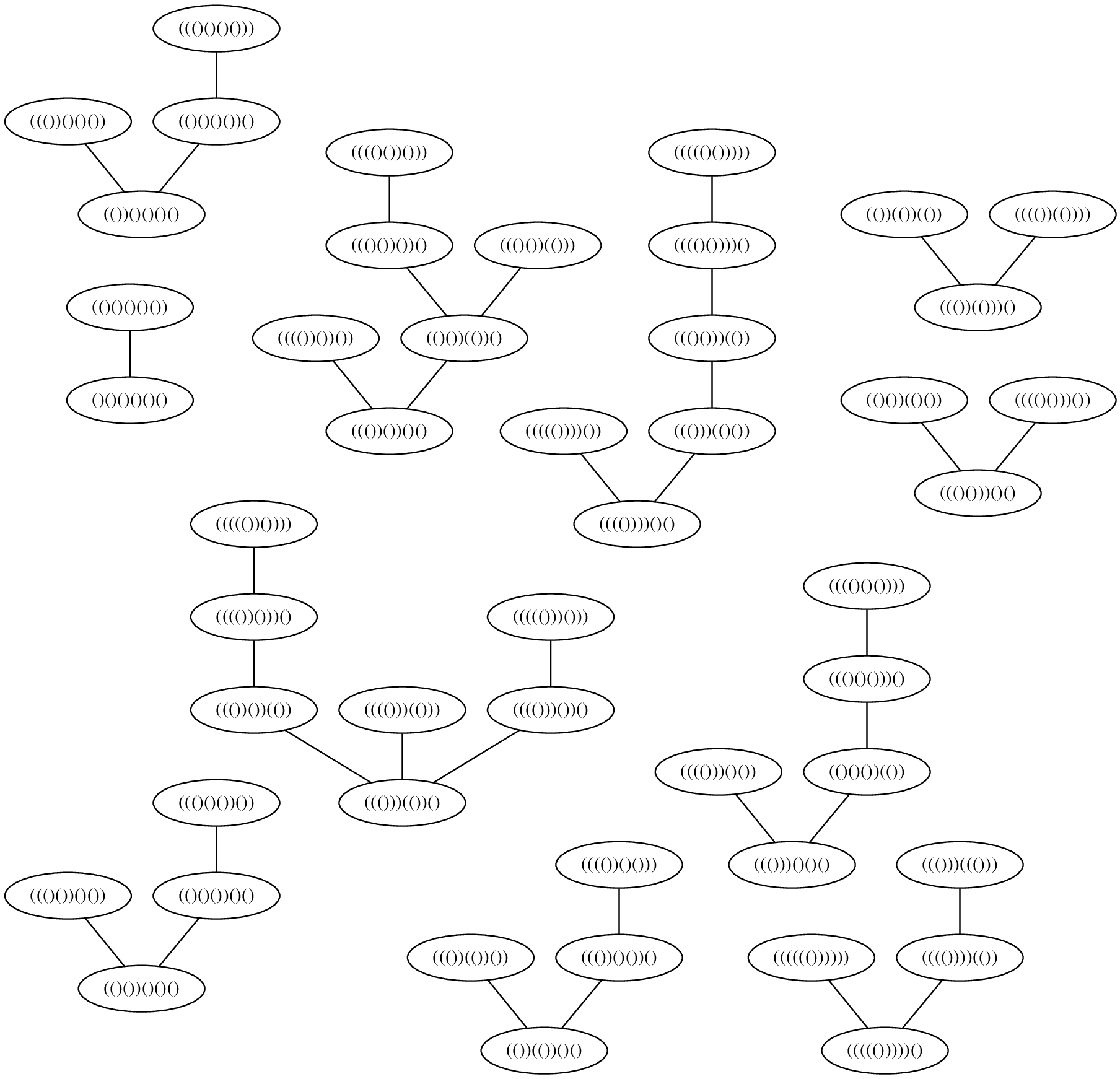}
\caption{The 11 clusters of grouping the $|\mathbb{C}_6|=48$
expressions with 6 pairs of parentheses into clusters.
}
\label{fig.C6}
\end{figure}

Any expression of nested circles in the plane can be interpreted
as a set of circles embedded in a sphere surface: draw a big
circumscribing circle around the set of all circles and interpret
it as an equator of the sphere. This defines a mapping of sets
of circle topologies of the plane onto one circle topology of the sphere, because
topologies that are related by flipping the interior and exterior
region of a circle are no longer distinct on the sphere. In our notation
of nested parentheses, such a flip starts with a nested expression \texttt{A ( B )},
where \texttt{A} and \texttt{B} are well-formed (potentially empty) subexpressions. After
embedding in the sphere, the closing parenthesis can be torn across the back surface of the
sphere to the opposite side, ending up with \texttt{ ( A ) B}. There are as many
flip operations as there are factors in the expression because one can move
any of them to the right before the flip---although some of their images
may be the same because factors may be equal, and although in some cases
the image may be the same as the original expression.
\begin{exa}
The flip operation on \texttt{()()()} gives \texttt{(()())}.
\end{exa}
\begin{exa}
The flip operation on \texttt{(())()} gives \texttt{((()))}
if the \texttt{()} factor is flipped. It gives \texttt{()(())}
when the \texttt{(())} factor is flipped; in that case the image is the
same as the original, because the order of the two factors does not matter here.
\end{exa}
Figures \ref{fig.C4}--\ref{fig.C6} illustrate for $N=4$--$6$ how expressions 
transform under the flip-transform: edges in the graphs mean that
the expression on one node is transformed to the expression of the other node
by a flip-transformation.


So the $\mathbb{C}_N$ circle sets in the plane can be sorted into clusters
which assemble all expressions that are mutually convertible by a chain
of flip transformations. The number of clusters in $\mathbb{C}_N$
equals the number of topologically distinct circle sets embedded in
the sphere surface.

\subsection{Sets of Nested Circles and Squares}
If the geometric figures have $k$ distinct hollow shapes---
for example circles and squares with $k=2$---the methods of 
circumscribing and placing side by side generalize the rules.
There are $k$ different ways of forming a single compound object
from a set of objects with one element less because there are
$k$ options for the outermost shape. An upper-left index $k$
specifies how many shapes are available.
(\ref{eq.CN1}) and (\ref{eq.CNf}) turn into
\begin{equation}
|^k\mathbb{C}_N^{(1)}|=k|^k\mathbb{C}_{N-1}|,
\label{eq.CNkof1}
\end{equation}
\begin{equation}
|^k\mathbb{C}_N^{(f)}| = \sum_{\pi(N): N=\{N_1^{c_1};N_2^{c_2};\ldots N_f^{c_f}\}} 
\prod_{j=1}^f
\binom{|^k\mathbb{C}_{N_j}^{(1)}|+c_j-1}{c_j}.
\end{equation}
The implicit equation of the generating function for the total
number of topologies in the plane is
\cite{LerouxASMQ16}
\begin{equation}
^kC(z)= \exp\left(k\sum_{j\ge 1} z^j\,\,{^kC(z^j)}/j\right).
\end{equation}

\begin{exa}
For $k=2$, $|^2\mathbb{C}_2|=7$ configurations exist:
a circle inside a circle, a circle inside a square, a square inside
a circle, a square inside a square, two disjoint circles, two disjoint squares,
or a separated square and circle.
\end{exa}
The case with $k=2$ shapes is further illustrated in Table \ref{tab.C2},
the case with $k=3$ shapes in Table \ref{tab.C3}.
The values on the diagonals are
\begin{equation}
|^k\mathbb{C}_N^{(N-1)}|=\binom{N+k-1}{k-1}.
\end{equation}
Row sums $|^k\mathbb{C}_N|$ are Euler transforms of the columns
$|^k\mathbb{C}_N^{(1)}|$.
\begin{table}
\begin{tabular}{r|r|rrrrrrrrrrrrrrrrrrrrr}
$N$ & $ |^2\mathbb{C}_N|$ & 1 & 2 & 3 & 4 & 5 & 6 & 7 &8 &9\\
\hline
1 & 2& 2  \\
2 & 7& 4 & 3 \\
3 & 26& 14 & 8 & 4 \\
4 & 107& 52 & 38 & 12 & 5 \\
5 & 458& 214 & 160 & 62 & 16 & 6 \\
6 & 2058& 916 & 741 & 288 & 86 & 20 & 7 \\
7 & 9498& 4116 & 3416 & 1408 & 416 & 110 & 24 & 8 \\
8 & 44947& 18996 & 16270 & 6856 & 2110 & 544 & 134 & 28 & 9 \\
9 & 216598& 89894 & 78408 & 34036 & 10576 & 2812 & 672 & 158 & 32 & 10 \\
\end{tabular}
\caption{The counts $|^2\mathbb{C}_N|$ and
$|^2\mathbb{C}_N^{(f)}|$ of nested nonintersecting circles and squares \cite[A000151,A038055,A271878]{EIS}.
}
\label{tab.C2}
\end{table}

\begin{table}
\begin{tabular}{r|r|rrrrrrrrrrrrrrrrrrrrr}
$N$ & $ |^3\mathbb{C}_N|$ & 1 & 2 & 3 & 4 & 5 & 6 & 7 & 8 & 9\\
\hline
1 &  3& 3 \\
2 &  15& 9 & 6 \\
3 &  82& 45 & 27 & 10 \\
4 &  495& 246 & 180 & 54 & 15 \\
5 &  3144& 1485 & 1143 & 405 & 90 & 21 \\
6 &  20875& 9432 & 7704 & 2856 & 720 & 135 & 28 \\
7 &  142773& 62625 & 52731 & 20682 & 5385 & 1125&  189&  36 \\
8 &  1000131& 428319 & 369969 & 150282 & 40914&  8730&  1620&  252&  45 \\
9 &  7136812& 3000393 & 2638332 & 1104702 & 309510&  68400&  12891&  2205&  324&  55 \\
\end{tabular}
\caption{The counts $|^3\mathbb{C}_N|$ and
$|^3\mathbb{C}_N^{(f)}|$ of nested nonintersecting circles, squares and triangles \cite[A006964,A038059,A271879]{EIS}.
}
\label{tab.C3}
\end{table}

\section{Topologically Distinct Circle Sets, One Circle Marked}\label{sec.M}
\subsection{Base-4 Notation}
In Section \ref{sec.C} the circles are moving without intersecting and 
qualitatively equal. We move on to the combinatorics of circle sets
where one of them is marked (for example by a unique color, by morphing
it into an ellipse or replacing it by a square). The equivalent modification
in the commutative algebra is to introduce another symbol, a pair of brackets \texttt{[]},
to locate the modified evaluation of a subexpression.
Factorization is defined as before, and the marked circle is not intersecting
with any of the other circles as before.

\begin{defn}
$\mathbb{M}_N^{(f)}$ is the set of $N$ circles with $f$ factors,
one of these circles marked.
\end{defn}

\begin{eqnarray}
\mathbb{M}_N &=& \bigcup_{f=1}^N \mathbb{M}_N^{(f)};\\
|\mathbb{M}_N| &=& \sum_{f=1}^N |\mathbb{M}_N^{(f)}|;\quad |\mathbb{M}_0|=1.
\end{eqnarray}

\begin{exa}\label{exa.M}
\begin{eqnarray}
\mathbb{M}_1 &=& \mathbb{M}_1^{(1)} = \{[]\}; \\
\mathbb{M}_2^{(1)} &=& \{[()], ([])\}; \\
\mathbb{M}_2^{(2)} &=& \{[]()\}; \\
\mathbb{M}_3^{(1)} &=& \{ (([])), ([()]), ([]()), [()()], [(())] \}; \\
\mathbb{M}_3^{(2)} &=& \{ ([])(), [()](), [](()) \}; \\
\mathbb{M}_3^{(3)} &=& \{ []()() \}; \\
\mathbb{M}_4^{(1)} &=& \{ 
((())[]),
((([]))),
(([()])),
(([])()),
(([]())),
([()()]), \nonumber \\
&&
([()]()),
([(())]),
([]()()),
[()()()],
[(())()],
[(()())],
[((()))]
\}; \\
\mathbb{M}_4^{(2)} &=& \{ 
(()())[],
((()))[],
(([]))(),
([()])(),
([])(()), \nonumber \\ &&
([]())(),
[()()](),
[()](()),
[(())]()
\};
\end{eqnarray}
\end{exa}

\subsection{Recurrences}
\begin{table}
\begin{tabular}{r|r|rrrrrrrrrrrrrrrrrrrrr}
$N$ & $ |\mathbb{M}_N|$ & 1 & 2 & 3 & 4 & 5 & 6 & 7 & 8 & 9 & 10  & 11\\
\hline
1 & 1&     1& \\
2 & 3&     2&     1&  \\
3 & 9&     5&     3&     1&  \\
4 & 26&    13&     9&     3&     1&  \\
5 & 75&    35&    26&    10&     3&     1&  \\
6 & 214&    95&    75&    30&    10&     3&     1&  \\
7 & 612&   262&   214&    91&    31 &   10&     3&     1&  \\
8 & 1747&   727&   612&   268&    95&    31&    10&     3&     1&  \\
9 & 4995&  2033&  1747&   790&   284&    96&    31&    10&     3&     1&  \\
10 & 14294&  5714&  4995&  2308&   848&   288&    96&    31&    10&     3&     1&  \\
11 & 40967& 16136& 14294&  6737&  2506&   864&   289&    96&    31&    10&     3&     1 
\end{tabular}
\caption{The number of nonintersecting circles with one of them marked.
$|\mathbb{M}_N|$ are the row sums and $|\mathbb{M}_N^{(f)}|$ the entries with $f$ factors
\cite[A000243,A000107]{EIS}.
}
\label{tab.M}
\end{table}

An overview of how many distinct arrangements exist is given in Table \ref{tab.M}.
The set of $\mathbb{M}_N^{(1)}$ is created by either (i) wrapping
an expression without a marked sphere into a bracket, or by (ii)
wrapping an expression that already contains a marked sphere into a pair
of parentheses:
\begin{equation}
|\mathbb{M}_N^{(1)}|
=
|\mathbb{C}_{N-1}|
+
|\mathbb{M}_{N-1}|.
\label{eq.Mrec1}
\end{equation}

For expressions with $f\ge 2$ factors we may always move the factor with the bracket
to some pivotal (say the leftmost) factor because the order of factors does
not matter in Sections \ref{sec.C} and \ref{sec.M}. That pivotal factor 
needs, say, $N'$ pairs of parentheses (including the marked),
and all the other factors may be varied as a set of the $\mathbb{C}$ type:
\begin{equation}
|\mathbb{M}_N^{(f)}|
=
\sum_{N'=1}^{N-1}
|\mathbb{M}_{N'}^{(1)}|\,
|\mathbb{C}_{N-N'}^{(f-1)}|,\quad f\ge 2
.
\label{eq.Mrec}
\end{equation}

\begin{defn}
The generating function of the topologies of non-intersecting
circles with one marked is
\begin{equation}
M(z)=\sum_{N\ge 0} |\mathbb{M}_N|z^N .
\label{eq.Mz}
\end{equation}
\end{defn}

If we sum on both sides of (\ref{eq.Mrec}) over $f$, insert
(\ref{eq.Cunion}) for the sum over the $\mathbb{C}$ and (\ref{eq.Mrec1})
to eliminate the $\mathbb{M}_{N'}^{(1)}$ on the right hand side,
Jovovic's relation shows up \cite[A000243]{EIS}
\begin{equation}
M(z)=1+\frac{zC^2(z)}{1-zC(z)}.
\end{equation}

For sufficiently large $N$ the count with $N-1$ factors is
$|\mathbb{M}_N^{(N-1)}|=3$ because the set contains the expressions
of the form $[()]()()\cdots$, $([])()()\cdots$,
and $[](())()()\cdots$.

\subsection{Inner Void Circle}
There is a subset of expressions $\mathbb{M}_N^{(v)}\subseteq \mathbb{M}_N$---indicated
with an upper $v$ like \emph{void}---where the bracket does not contain
any subexpression with parentheses, i.e., where the marked circle does not
circumscribe any other circle.

\begin{defn}
$\mathbb{M}_N^{(f,v)}$ is the set of $N$ circles with $f$ factors,
where one of these circles is marked and does not contain other circles.
\end{defn}

\begin{eqnarray}
\mathbb{M}_N^{(v)} &=& \bigcup_{f=1}^N \mathbb{M}_N^{(f,v)} \subseteq \mathbb{M}_N^{(f)};\\
|\mathbb{M}_N^{(v)}| &=& \sum_{f=1}^N |\mathbb{M}_N^{(f,v)}|.
\end{eqnarray}
\begin{exa}
Removing in Example \ref{exa.M} the expressions where the bracket pair
embraces other parentheses yields:
\begin{eqnarray}
\mathbb{M}_1^{(v)} &=& \mathbb{M}_1^{(1,v)} = \{[]\}; \\
\mathbb{M}_2^{(1,v)} &=& \{([])\}; \\
\mathbb{M}_2^{(2,v)} &=& \{[]()\}; \\
\mathbb{M}_3^{(1,v)} &=& \{ (([])), ([]()), \}; \\
\mathbb{M}_3^{(2,v)} &=& \{ ([])(), [](()) \}; \\
\mathbb{M}_3^{(3,v)} &=& \{ []()() \}; \\
\mathbb{M}_4^{(1,v)} &=& \{ 
((())[]),
((([]))),
(([])()),
(([]())),
([]()()),
\}; \\
\mathbb{M}_4^{(2,v)} &=& \{ 
(()())[],
((()))[],
(([]))(),
([])(()),
([]())()
\}.
\end{eqnarray}
\end{exa}

\begin{table}
\begin{tabular}{r|r|rrrrrrrrrrrrrrrrrrrrr}
$N$ & $ |\mathbb{M}_N^{(v)}|$ & 1 & 2 & 3 & 4 & 5 & 6 & 7 & 8 & 9 & 10  & 11\\
\hline
1 & 1&     1 \\
2 & 2&     1&     1 \\
3 & 5&     2 &    2  &   1 \\
4 & 13&     5 &    5  &   2 &    1 \\
5 & 35&    13 &   13  &   6  &   2 &    1 \\
6 & 95&    35 &   35  &  16  &   6 &    2  &   1 \\
7 & 262&    95 &   95 &   46  &  17&     6 &    2 &    1 \\
8 & 727&   262  & 262  & 128   & 49 &   17 &    6 &    2  &   1 \\
9 & 2033&   727 &  727  & 364  & 139&    50 &   17 &    6 &    2 &    1 \\
10 & 5714&  2033 & 2033 & 1029  & 401 &  142&    50 &   17 &    6 &    2 &    1 \\
11 & 16136&  5714 & 5714 & 2930  &1147&   412&   143&    50&    17 &    6 &    2 &    1 \\
12 & 45733& 16136 &16136  &8344  &3299&  1184&   415 &  143&    50 &   17 &    6 &    2 &    1 \\
\end{tabular}
\caption{The number of nonintersecting circles, one marked.
$|\mathbb{M}_N^{(v)}|$
and
$|\mathbb{M}_N^{(f,v)}|$ for $N,f \ge 1$
\cite[A000107]{EIS}.
}
\label{tab.Mv}
\end{table}

The topologies with that scenario are counted in Table \ref{tab.Mv}.
With the same argument as in Equation (\ref{eq.Mrec}), scenarios
with an empty bracket need to locate the bracket at some fixed factor, and let the other factors generate
all possible diagrams with the remaining parentheses:
\begin{equation}
|\mathbb{M}_N^{(f,v)}|
=
\sum_{N'=1}^{N-1}
|\mathbb{M}_{N'}^{(1,v)}|\,
|\mathbb{C}_{N-N'}^{(f-1)}|
.
\label{eq.Mrecv}
\end{equation}

On the diagonals of Tables \ref{tab.M} and \ref{tab.Mv} we find
\begin{equation}
|\mathbb{M}_N^{(N)}|
=
|\mathbb{M}_N^{(N,v)}|
=1,
\end{equation}
because the only expressions with as many factors as circles
is the product of singletons, $\mathbb{M}_N^{(N)} = \mathbb{M}_N^{(N,v)} = \{[]()()\cdots ()\}$.

The number in column $\mathbb{M}_N^{(1)}$ in Table \ref{tab.Mv}
duplicates the total of the previous row:
\begin{equation}
|\mathbb{M}_N^{(1,v)}|=|\mathbb{M}_{N-1}^{(v)}|
.
\end{equation}
This is easily understood because each element of the set $\mathbb{M}_N^{(1)v}$
is created by surrounding the expression of an element of the set $\mathbb{M}_{N-1}^{(v)}$
by a pair of (non-marked) parentheses, so the ``void'' within the bracket is conserved.

In a similar manner $|\mathbb{M}_N^{(1,v)}| =
|\mathbb{M}_N^{(2,v)}|$ is understood by ``peeling off'' the outermost pair of parentheses
of the element of $\mathbb{M}_N^{(1,v)}$
and placing it as an extra factor $()$ aside from the peeled expression. This association
works because the outermost pair of parentheses is never the bracket.

In summary, all entries of Table \ref{tab.M} and \ref{tab.Mv} can be recursively
generated from Table \ref{tab.C} with the aforementioned 4 formulas.

\begin{rem}
The serialized representation of the circle sets with two types of parentheses
on a computer is possible by moving from the binary digit representation
of Sections \ref{sec.P} and \ref{sec.C} to a base-4 representation 
$)\mapsto 0$,
$(\mapsto 1$,
$]\mapsto 2$,
$[\mapsto 3$. The mapping is $[[]]\mapsto 3322_4$, $([[()]]())\mapsto 1331022100_4$, for example.
\end{rem}

\section{Circle Sets With One Pair intersecting} 
\subsection{Serialized Notation}\label{sec.Xser}
Another derivative of the non-intersecting circle sets of Section \ref{sec.C}
are circle sets where exactly one pair of circles intersects at two points
of their rims. 

These two intersecting circles are a natural reference frame for the other $N-2$.
In the serialized notation we introduce the expression $[[]]$ with
two bracket pairs to indicate crossing of the rims of the first, then of the
second circle, then leaving the first and finally leaving the second.
The notation provides 5 regions that host the $N-2$ remaining circles. The
well-formed general expression 
is \textit{reg4}$[$\textit{reg3}$[$\textit{reg2}$]$\textit{reg1}$]$\textit{reg0}
if the regions are enumerated 0--4.

$\textit{reg4}[\textit{reg3}[\textit{reg2}]\textit{reg1}]\textit{reg0}\mapsto$
\begin{picture}(6,1.5)
\put(0.5,0.8){reg4}
\put(2.2,0.8){\circle{1.6}}
\put(1.7,1.0){reg}
\put(1.7,0.7){3}
\put(2.3,0.8){reg}
\put(2.4,0.5){2}
\put(2.9,0.8){\circle{1.6}}
\put(2.9,0.9){reg}
\put(3.0,0.6){1}
\put(4.5,0.8){reg0}
\end{picture}

The serialized notation is well-suited for computerized managing,
but again has the drawback that the freedom of moving circle sets
around as long as no new intersections are induced is not strictly enforced.
We add the following constraints to the serialized notation to avoid over-counting
those circle sets with two intersections:
\begin{enumerate}
\item
The regions \textit{reg1}, \textit{reg2} and \textit{reg3} host members
of the $\mathbb{C}_N$ collection. This basically ensures that their
circle sets do not introduce intersections by peeking beyond the enclosures
defined by the bracket pair.
Note that no such rule is enforced on \textit{reg0} and \textit{reg4}
because we allow the crossing circles to be inside other circles; so
an expression like $([[]])$ is well-formed, although the isolated left
and right parentheses are not individually members of $\mathbb{P}$.
\item
If the entire core region of the crossing circles is removed---leaving
the concatenated expression \textit{reg4}\textit{reg0}---this must be a well-formed $\mathbb{P}$
expression. This ensures that circles that rotate in the space outside
the crossing circles are considered equivalent; eventually
expressions like $(()[[]])$ and $([[]]())$ for example are counted only once.
\item
From the two expressions obtained by swapping \textit{reg1} and \textit{reg3}
only one is admitted. These are the regions inside one of the intersecting
circles but not in the intersection. The rule ensures that a sort
of mirror operation at the center of the intersection---which does not
change the topology---is admitted only once in the circle sets.
\end{enumerate}

\begin{defn}
$\mathbb{X}_N^{(f)}$ is the set of $N$ circles with $f$ factors,
two circle rims intersecting in two points.
\end{defn}

\begin{eqnarray}
\mathbb{X}_N &=& \bigcup_{f=1}^N \mathbb{X}_N^{(f)};\\
|\mathbb{X}_N| &=& \sum_{f=1}^N |\mathbb{X}_N^{(f)}|; \quad |\mathbb{X}_1|=0.
\end{eqnarray}

\begin{exa}
$([[]])() \mapsto$
\begin{picture}(6,1.5)
\put(1.2,0.5){\circle{0.5}}
\put(1.5,0.5){\circle{0.5}}
\put(1.35,0.5){\circle{3.0}}
\put(2.5,0.5){\circle{0.5}}
\end{picture}
$([[]()]) \mapsto$
\begin{picture}(6,1.5)
\put(1.1,0.5){\circle{0.7}}
\put(1.6,0.5){\circle{0.7}}
\put(1.35,0.5){\circle{3.0}}
\put(1.7,0.5){\circle{0.3}}
\end{picture}
$[[]](()) \mapsto$
\begin{picture}(6,1.5)
\put(1.1,0.5){\circle{0.7}}
\put(1.6,0.5){\circle{0.7}}
\put(2.6,0.5){\circle{1.0}}
\put(2.6,0.5){\circle{0.5}}
\end{picture}
$[[()]()]() \mapsto$
\begin{picture}(6,1.5)
\put(1.1,0.5){\circle{2.0}}
\put(1.8,0.5){\circle{2.0}}
\put(2.9,0.5){\circle{0.5}}
\put(2.1,0.5){\circle{0.5}}
\put(1.45,0.5){\circle{0.4}}
\end{picture}
\end{exa}

\begin{exa}\label{ex.X}
\begin{eqnarray}
\mathbb{X}_2&=& \mathbb{X}_2^{(1)} = \{[[]]\}; \\
\mathbb{X}_2^{(2)} &=& \{\};\\
\mathbb{X}_3^{(1)} &=&
\{([[]]),
[[()]],
[[]()]\}; \\
\mathbb{X}_3^{(2)} &=&
\{[[]]()\}; \\
\mathbb{X}_3^{(3)} &=& \{\};\\
\mathbb{X}_4^{(1)} &=&
\{(([[]])),
([[()]]),
([[]()]),
([[]]()),
\nonumber \\ &&
[()[]()],
[[()()]],
[[()]()],
[[(())]],
[[]()()],
[[](())]\}; \\
\mathbb{X}_4^{(2)} &=&
\{([[]])(),
[[()]](),
[[]()](),
[[]](())\}; \\
\mathbb{X}_4^{(3)} &=&
\{[[]]()()\}; \\
\mathbb{X}_4^{(4)} &=& \{\};
\end{eqnarray}
\end{exa}
\subsection{Recurrences}\label{sec.Xnrecur}
Table \ref{tab.X} shows how many expressions are in the sets $\mathbb{X}_N$
and $\mathbb{X}_N^{(f)}$. The first three values of $|\mathbb{X}_N|$ are mentioned
in the Encyclopedia of Integer Sequences \cite[A261070]{EIS}.
\begin{table}
\begin{tabular}{r|r|rrrrrrrrrrrrrr}
$N$ &  $|\mathbb{X}_N|$  & 1 & 2 & 3 & 4 & 5 & 6 & 7 & 8 & 9 & 10\\
\hline
2 & 1 &1 &0 &\\ 
3 & 4 &3 &1 &0 &\\ 
4 & 15 &10 &4 &1 &0 &\\ 
5 & 50 &30 &15 &4 &1 &0 &\\ 
6 & 162 &91 &50 &16 &4 &1 &0 &\\ 
7 & 506 &268 &162 &55 &16 &4 &1 &0 &\\ 
8 & 1558 &790 &506 &185 &56 &16 &4 &1 &0 &\\ 
9 & 4727 &2308 &1558 &594 &190 &56 &16 &4 &1 &0 &\\ 
10 & 14227&   6737&  4727&  1878&   617&   191&    56&    16&     4&     1&0 \\ 
11 & 42521&  19609& 14227&  5825&  1970&   622&   191&    56&    16&     4 &1  &0 
\end{tabular}
\caption{Topologically distinct sets of $N$ circles with one pair intersecting, total (row sums) $|\mathbb{X}_N|$
and $|\mathbb{X}_N^{(f)}|$
classified according to number of factors $1\le f\le N$.}
\label{tab.X}
\end{table}

Obviously $|\mathbb{X}_N^{(N)}|=0$ and $|\mathbb{X}_N^{(N-1)}|=1$ because
we always spend two circles in the bracket---which does not factorize---and
the expression $[[]]()()()\cdots$ is the only member of $\mathbb{X}_N^{(N-1)}$.

For sufficiently large $N$ there are $|\mathbb{X}_N^{(N-2)}|=4$ expressions,
namely
$[[()]]()()\cdots$,
$[[]()]()()\cdots$,
$[[]](())()\cdots$, and
$([[]])()()\cdots$ with $N-3$ trailing isolated circles.

The argument of isolating the factor that contains the bracket pair
that led to Equation (\ref{eq.Mrec}) remains valid, so
\begin{equation}
|\mathbb{X}_N^{(f)}|
=
\sum_{N'=1}^{N-1}
|\mathbb{X}_{N'}^{(1)}|\,
|\mathbb{C}_{N-N'}^{(f-1)}|,\quad f\ge 2
.
\label{eq.Xoff}
\end{equation}

The dismantling of the sole factor of an expression of $\mathbb{X}_N^{(1)}$
that contains the two brackets shows two variants:
if the outer parentheses are the round parenthesis, the expression has been
formed by embracing any expression with $N-1$ circles, which contributes $|\mathbb{X}_{N-1}|$.
If alternatively the expression is of the form stripped down to where \textit{reg4}
and \textit{reg0} are empty, we count the number of ways of construction \textit{reg3},
\textit{reg2} and \textit{reg1} with a total of $N-2$ circles by a function $D_{n-2}$:
\begin{equation}
|\mathbb{X}_N^{(1)}| = |\mathbb{X}_{N-1}| + D_{N-2}.
\label{eq.XofD}
\end{equation}
\begin{equation}
D_N  = 1, 2, 6, 15, 41, 106, 284, 750, 2010, 5382, 14523, 39290\ldots;\quad N\ge 0.
\end{equation}
\begin{exa}
The 6 expressions that contribute to $D_2=6$ are
$ [[]()()]$,
$[[](())]$,
$[[()]()]$,
$[()[]()]$,
$[[(())]]$, and
$[[()()]]$.
\end{exa}

The distribution of the $N$ circles over \textit{reg3}, \textit{reg2} and \textit{reg1}
has no further restrictions to place any member of $\mathbb{C}$
into \textit{reg2}, which reduces $D$ by composition to another function $\hat D$
of the form
\begin{equation}
D_N=\sum_{N'=0}^{N} |\mathbb{C}_{N'}| \hat D_{N-N'}.
\label{eq.DNsplit}
\end{equation}
$\hat D_N$ counts the number of ways of placing $N$ circles 
in total into \textit{reg3}
and \textit{reg1} such that each expression is a member of $\mathbb{C}$ and such
that the third rule of Section \ref{sec.Xser} of counting  only the ``ordered'' pairs
is obeyed. If $N$ is odd, the expressions in two regions necessarily differ because
they must have a different number of circles, so the rule may for example be
implemented by putting always the expression with the lower number into one region:
\begin{equation}
\hat D_{N,odd} = \sum_{N'=0}^{\lfloor N/2\rfloor} |\mathbb{C}_{N'}|\,|\mathbb{C}_{N-N'}|.
\end{equation}
If $N$ is even, an additional format appears where the expressions in \textit{reg3}
and \textit{reg1} have the same number of circles. Because these elements of $|\mathbb{C}_{N/2}|$
may be put into a strict order, the triangular number with that argument counts
the ``non-ordered'' pairs of these:
\begin{equation}
\hat D_{N,even} = \sum_{N'=0}^{N/2} |\mathbb{C}_{N'}|\,|\mathbb{C}_{N-N'}|
+ \frac{|\mathbb{C}_{N/2}|(|\mathbb{C}_{N/2}|+1)}{2}.
\end{equation}
In terms of the generating functions (\ref{eq.Cz}), (\ref{eq.Mz}) and
\begin{equation}
\hat D(z) = \sum_{N\ge 0} \hat D_N z^N,\quad
D(z) = \sum_{N\ge 0} D_N z^N,\quad
\end{equation}
this type of half convolution in the previous two equations
may be summarized as \cite[A027852]{EIS}
\begin{equation}
\hat D(z)=\frac12[C(z)^2+C(z^2)].
\label{eq.Dogf}
\end{equation}
\begin{rem}
The symmetry enforced to the contents of \textit{reg1} and \textit{reg3} is the symmetry
of the cyclic group of order 2\@. The cycle index for this group is $(t_1^2+t_2)/2$ \cite[I60]{Flajolet}. Substitution
of $t_j\mapsto C(z^j)$ gives the same result \cite[p. 252]{Comtet}.
\label{rem.cycC2}
\end{rem}
\begin{equation}
\hat D_N = 1, 1, 3, 6, 16, 37, 96, 239, 622, 1607, 4235, 11185, 29862\ldots ;\quad N\ge 0.
\end{equation}
The convolution (\ref{eq.DNsplit}) turns into a product of the
generating functions:
\begin{equation}
D(z)= C(z)\hat D(z).
\end{equation}
\begin{exa}
$\hat D_1=1$ is the size of the set $\{[[]()]\}$.
\end{exa}
\begin{exa}
$\hat D_2=3$ is the size of the set $\{[()[]()], [[]()()],[[](())]\}$.
\label{exa.hatD2}
\end{exa}
\begin{exa}
$\hat D_3=6$ is the size of the set $\{
[[]()()()],
[[](())()],$
$[[]((()))],$ 
$[[](()())],$ 
$[()[]()()],$ 
$[()[](())]
\}$.
\label{exa.hatD3}
\end{exa}

Summing (\ref{eq.Xoff}) over $f$ and using (\ref{eq.XofD}) leads to
\begin{equation}
X(z) =1+\frac{z^2D(z)C(z)}{1-zC(z)}.
\label{eq.Xogf}
\end{equation}

\subsection{Pair of Touching Circles} 
If there are no further circles in the area of the intersection, the
two intersecting circles may be moved apart until they touch in a single point.
These borderline cases are destilled from the previous analysis by counting
expressions only where the two inner brackets appear side by side, i.e,
where \textit{reg2} is empty.
We call these sets of configurations $\mathbb{X}_N^{(f,t)}$ where the label $t$ indicates
touching.

\begin{defn}
$\mathbb{X}_N^{(f,t)}$ is the set of $N$ circles with $f$ factors,
two circle rims touching at one point.
\end{defn}

\begin{eqnarray}
\mathbb{X}_N^{(t)} &=& \bigcup_{f=1}^N \mathbb{X}_N^{(f,t)}\subseteq \mathbb{X}_N;\\
|\mathbb{X}_N^{(t)}| &=& \sum_{f=1}^N |\mathbb{X}_N^{(f,t)}| \le |\mathbb{X}_N|;
\end{eqnarray}

\begin{exa}
If we remove the expressions from Example \ref{ex.X} where other circles
appear within the innermost of the two square brackets, the following list emerges:
\begin{eqnarray}
\mathbb{X}_2^{(t)}&=& \mathbb{X}_2^{(1,t)} = \{[[]]\}; \\
\mathbb{X}_2^{(2,t)} &=& \{\};\\
\mathbb{X}_3^{(1,t)} &=&
\{([[]]),
[[]()]\}; \\
\mathbb{X}_3^{(2,t)} &=&
\{[[]]()\}; \\
\mathbb{X}_3^{(3,t)} &=& \{\};\\
\mathbb{X}_4^{(1,t)} &=&
\{(([[]])),
([[]()]),
([[]]()),
\{[()[]()],
[[]()()],
[[](())]\}; \\
\mathbb{X}_4^{(2,t)} &=&
\{([[]])(),
[[]()](),
[[]](())\}; \\
\mathbb{X}_4^{(3,t)} &=&
\{[[]]()()\}; \\
\mathbb{X}_4^{(4)} &=& \{\};
\end{eqnarray}
\end{exa}

Table \ref{tab.Xt} shows how many expressions are in the sets $\mathbb{X}_N^{(t)}$
and $\mathbb{X}_N^{(f)t}$.
\begin{table}
\begin{tabular}{r|r|rrrrrrrrrrrrrr}
$N$ &  $|\mathbb{X}_N^{(t)}|$  & 1 & 2 & 3 & 4 & 5 & 6 & 7 & 8 & 9 & 10\\
\hline
1 & 0 &      0 &\\
2 & 1 &      1 &     0 &\\
3 & 3 &      2 &     1 &     0 &\\
4 & 10 &      6 &     3 &     1 &     0 &\\
5 & 30 &     16 &    10 &     3 &     1 &     0 &\\
6 & 91 &     46 &    30 &    11 &     3 &     1 &     0 &\\
7 & 268 &    128 &    91 &    34 &    11 &     3 &     1 &     0 &\\
8 & 790 &    364 &   268 &   108 &    35 &    11 &     3 &     1 &     0 &\\
9 & 2308 &   1029 &   790 &   327 &   112 &    35 &    11 &     3 &     1 &     0 &\\
10 & 6737 &   2930 &  2308 &   992 &   344 &   113 &    35 &    11 &     3 &     1 &     0 &\\
11 & 19609 &   8344 &  6737 &  2962 &  1055 &   348 &   113 &    35 &    11 &     3 &     1 &     0 &\\
\end{tabular}
\caption{Topologically distinct sets of $N$ circles with one pair touching, total
and classified according to number of factors $1\le f\le N$ \cite[A269800]{EIS}.}
\label{tab.Xt}
\end{table}

As before
\begin{equation}
|\mathbb{X}_N^{(N,t)}|=0;\quad |\mathbb{X}_N^{(N-1,t)}|=1.
\end{equation}
The strategy of isolating the factor with the brackets that lead to (\ref{eq.XofD}) remains valid:
\begin{equation}
|\mathbb{X}_N^{(f,t)}|
=
\sum_{N'=1}^{N-1}
|\mathbb{X}_{N'}^{(1,t)}|\,
|\mathbb{C}_{N-N'}^{(f-1)}|,\quad f\ge 2
.
\end{equation}
The formula that distributed the $N-2$ circles within the three regions in the 
intersecting circles now needs to skip the cases where some of them
are in \textit{reg2}. And instead of (\ref{eq.XofD}) we immediately skip to
\begin{equation}
|\mathbb{X}_N^{(1,t)}| = |\mathbb{X}_{N-1}^{(t)}| + \hat D_{N-2}
\label{eq.XNof1t}
\end{equation}
and replace (\ref{eq.Xogf}) by the generating function
\begin{equation}
X^{(t)}(z) = \sum_{N\ge 0}X^{(t)}_N z^N=1+\frac{z^2\hat D(z)C(z)}{1-zC(z)}.
\end{equation}

\subsection{One or More Intersecting Pairs}\label{sec.tree}
The topologies of the members of the
sets $\mathbb{X}_N$
are mapped
onto rooted trees representing
the dependence of ``being a circle inside another''
as ``being a branch of a node that represents the
enclosing circle.'' The plane is the root of the
tree. There is no limit of how many branches a node
can have. Moving around circle clusters freely means
that the nodes are counted without a notion of order.
The sole exception---which distinguishes $\mathbb{C}_N$
from $\mathbb{X}_N$---is that the tree must have
a single node representing the intersecting circle pair
which has up to three nodes (the three regions) that partially
respect order because the circle clusters in \textit{reg2}
are topologically considered different from circle clusters
in \textit{reg1} or \textit{reg3}.

If one chops the node representing the plane off the tree,
it becomes a rooted forest, where the number of rooted trees
is the factor $f$ of the interpretation as nested parentheses.

The natural extension of these rules is to 
symmetrize the rules for the branches in that rooted forest, i.e.,
to allow circles and the three regions in the intersecting
circle pair to host any number of intersecting circle
clusters or intersecting circle pairs. The restriction
that remains is that intersection of more than two
circles are still not considered.

\begin{defn}
$^2\mathbb{X}_N$ is the set of topologies of $N$
nested circles in the plane where each circle intersects
with at most one other circle. $^2\mathbb{X}_N^{(f)}$
is the set of topologies of $N$ nested circles with $f$
factors, i.e., with $f$ of these objects that are
not inside any other of these objects.
\end{defn}
\begin{defn}
The generating function is
\begin{equation}
^2X(z) = \sum_{N\ge 0} |^2\mathbb{X}_N| z^N.
\end{equation}
\end{defn}

\begin{exa}
This is a circle bundle in $^2\mathbb{X}_7^{(1)}$ which is not
in $\mathbb{X}_7$:
\begin{picture}(6,1.5)
\put(1.0,0.9){\circle{0.5}}
\put(1.6,0.9){\circle{0.5}}
\put(1.35,0.9){\circle{3.0}}
\put(2.5,0.7){\circle{0.5}}
\put(2.6,1.0){\circle{0.5}}
\put(3.0,0.7){\circle{0.3}}
\put(2.6,0.9){\circle{3.0}}
\end{picture}
The outer pair of 2 intersecting circles contains another
pair of 2 intersecting circles (amongst others)
in one of its three regions.
\end{exa}

The grand book-keeping of placing these objects side
by side works as before, and the empty plane
is the unique way of not having any circles:
\begin{equation}
^2\mathbb{X}_N = \bigcup_{f=1}^N{} ^2\mathbb{X}_N^{(f)};
\quad
|^2\mathbb{X}_N| = \sum_{f=1}^N |^2\mathbb{X}_N^{(f)}|;
\quad
|^2\mathbb{X}_0| = 1.
\end{equation}
\begin{equation}
|^2\mathbb{X}_N^{(f)}| = \sum_{\pi(N): N=\{N_1^{c_1};N_2^{c_2};\ldots N_f^{c_f}\}} 
\prod_{j=1}^f
\binom{|^2\mathbb{X}_{N_j}^{(1)}|+c_j-1}{c_j};\quad f\ge 2.
\end{equation}
The difference starts where the objects at $f=1$
are dismantled. These are not the two types considered in
(\ref{eq.Mrec1}),
(\ref{eq.XofD})
or (\ref{eq.XNof1t})
nor the $k$ types as in (\ref{eq.CNkof1}).
The compound object is either a circle that hosts
the same type of objects with one circle less, or a pair of
intersecting circles with other objects of the
same type in their regions:
\begin{equation}
|^2\mathbb{X}_N^{(1)}| 
= |^2\mathbb{X}_{N-1}| + \bar D_{N-2}.
\label{eq.X2of1}
\end{equation}
$\bar D_{N-2}$ is the number of ways of
distributing objects of the $^2\mathbb{X}$ type
with a total of $N-2$ circles into three regions with
the symmetry rules of Section \ref{sec.Xnrecur}.

With the splitting rule of Section \ref{sec.Xnrecur}
the overlapping \textit{reg2} may contain any number
of the elements of $^2\mathbb{X}$ and the
other two regions share the remaining number of circles
as if the set was ordered. Copying from (\ref{eq.DNsplit}) and (\ref{eq.Dogf}),
\begin{equation}
\bar D_N = \sum_{N'=0}^N |^2\mathbb{X}_{N'}| \, \tilde D_{N-N'};
\label{eq.Dbardef}
\end{equation}
\begin{equation}
\tilde D(z) = \sum_{N\ge 0} \tilde D_N z^N 
= \frac12[^2X(z)^2+{}^2X(z^2)].
\end{equation}

\begin{table}
\begin{tabular}{r|r|rrrrrrrrrrrrrr}
$N$ & $ |^2\mathbb{X}_N|$ & 1 & 2 & 3 & 4 & 5 & 6 & 7 & 8 & 9 &\\
\hline
1 & 1&      1 & \\
2 & 3&      2 &      1 & \\
3 & 8&      5 &      2 &      1 & \\
4 & 27&     16 &      8 &      2 &      1 & \\
5 & 90&     53 &     26 &      8 &      2 &      1 & \\
6 & 330&    189 &    100 &     30 &      8 &      2 &      1 & \\
7 & 1225&    694 &    375 &    115 &     30 &      8 &      2 &      1 & \\
8 & 4729&   2642 &   1473 &    453 &    120 &     30 &      8 &      2 &      1 & \\
9 & 18554&  10270 &   5823 &   1827 &    473 &    120 &     30 &      8 &      2 &      1 & \\
10 & 74234&  40747 &  23479 &   7432 &   1936 &    479 &    120 &     30 &      8 &      2 &      1 & \\
11 & 300828& 164033 &  95618 &  30622 &   7954 &   1961 &    479 &    120 &     30 &      8 &      2 &      1 & \\
\end{tabular}
\caption{The number of topologies of nested $N$ circles intersecting at most as binaries, $|^2\mathbb{X}_N|$,
and the subcounts with $f$ factors, $|^2\mathbb{X}_N^{(f)}|$, $1\le f\le N$.
The row sums are the Euler transform of the column $f=1$.
}
\label{tab.X2}
\end{table}

\begin{equation}
\bar D_N = 
1, 2, 8, 26, 99, 364, 1417, 5541, 22193, 89799, 368160, 1523020,
\ldots, N\ge 0.
\end{equation}

The number of ways of distributing $N$ circles over
\textit{reg1} and \textit{reg3} of two intersecting circles is
\begin{equation}
\tilde D_N = 
1, 1, 4, 11, 41, 141, 537, 2041, 8042, 32028, 129780, 531331, 2198502,
\ldots N\ge 0.
\end{equation}
\begin{exa}
$\tilde D_2=4$ counts the three ways of Example \ref{exa.hatD2} plus
the one way of putting two intersecting circles in
one of the two regions.
\end{exa}
\begin{exa}
$\tilde D_3=11$ counts the six ways of Example \ref{exa.hatD3}
plus the following five ways that are new in $^2\mathbb{X}_N$:
\begin{enumerate}
\item
putting $[[]]()$ in one region,
\item
putting $[[]]$ in one region and $()$ in the other,
\item
putting $([[]])$ in one region,
\item
putting $[[]()]$ in one region,
\item
putting $[[()]]$ in one region.
\end{enumerate}
\end{exa}

\section{Outlook: 3-circle Intersections}
\subsection{6 New Topologies}
Adding 3-circle intersections introduces 6 new topologies 
beyond those of $^2\mathbb{X}_N$ \cite[A250001]{EIS}:
\begin{enumerate}
\item
The ``RGB spot diagram'' $^{3,1}\mathbb{X}_3$:
\begin{picture}(4,2.5)
\put(1.0,0.5){\circle{2.0}}
\put(2.0,0.5){\circle{2.0}}
\put(1.5,1.5){\circle{2.0}}
\put(2.0,0.5){1}
\put(1.5,1.5){2}
\put(1.0,0.5){3}
\end{picture}
The 7 regions inside the circles may be labeled by the circles that cover them: \textit{1}, \textit{12},
\textit{2}, \textit{23}, \textit{3}, \textit{13}, \textit{123}.
The symmetry of the diagram is established by three mirror lines that pass
through \textit{12}, \textit{23} and \textit{13}, and a symmetry for 
rotations by multiples
of 120$^\circ$ around the center. The symmetry group is the noncyclic group
of order 6\@. A permutation representation is (1)(23)
for the first generator and (123) for the second \cite{MatharVixra1504}.
The elements are
\begin{itemize}
\item
the unit element (1)(2)(3) which contributes 
$t_1^3$ to the cycle polynomial,
\item
the first generator which contributes
$t_1t_2$,
\item
the second generator which contributes $t_3$,
\item
the square of the second generator, (132), which contributes $t_3$,
\item
the element (12) which contributes $t_1t_2$, and
\item
the element (13)
which contributes $t_1t_2$.
\end{itemize}
The cycle index is $(t_1^3+3t_1t_2+2t_3)/6$.

\item
The torn version of this with an uncovered central area $^{3,2}\mathbb{X}_3$:
\begin{picture}(4,2.5)
\put(1.0,0.5){\circle{1.2}}
\put(2.0,0.5){\circle{1.2}}
\put(1.5,1.6){\circle{1.2}}
\put(2.0,0.5){1}
\put(1.5,1.5){2}
\put(1.0,0.5){3}
\end{picture}
The 7 regions inside the circles may be labeled by the circles that cover them: \textit{1}, \textit{12},
\textit{2}, \textit{23}, \textit{3}, \textit{13}, \textit{$\notin{\textit{123}}$}.
The symmetry is the same as for $^{3,1}\mathbb{X}_3$ above.
\item
A linear chain $^{3,3}\mathbb{X}_3$:
\begin{picture}(4,1.5)
\put(1.0,0.5){\circle{3.0}}
\put(1.0,0.5){3}
\put(2.0,0.5){\circle{3.0}}
\put(2.0,0.5){2}
\put(3.,0.5){\circle{3.0}}
\put(3.0,0.5){1}
\end{picture}
The 5 regions inside the circles may be labeled by the circles that cover them: \textit{1}, \textit{12},
\textit{2}, \textit{23}, \textit{3}.
The symmetry is the same left-right mirror symmetry as in Remark \ref{rem.cycC2};
the cycle index is $(t_1^2+t_2)/2$.
\item
The left-right compressed version of the previous diagram, $^{3,4}\mathbb{X}_3$:
\begin{picture}(5,1.5)
\put(0.9,0.5){\circle{3.0}}
\put(1.5,0.5){\circle{3.0}}
\put(2.1,0.5){\circle{3.0}}
\put(2.3,0.5){1}
\put(1.4,0.9){2}
\put(0.5,0.5){3}
\end{picture}
The 7 regions inside the circles may be labeled by the circles that cover them: \textit{1}, \textit{12},
\textit{123}, \textit{23}, \textit{3}, \textit{$\overline{\textit{2}}$}, \textit{$\underline{\textit{2}}$},
using overline and underline to register the upper and lower regions
of the pieces of circle 2.
The appearance of the regions 
$\overline{\textit{2}}$
and $\underline{\textit{2}}$ introduces an additional up-down mirror symmetry.
The symmetry group is the noncyclic group of order 4, which
has the generators (34) and (12) \cite{MatharVixra1504}.
The elements are
\begin{itemize}
\item
the unit element (1)(2)(3)(4) which contributes 
$t_1^4$ to the cycle polynomial,
\item
the first generator which contributes
$t_1^2t_2$,
\item
the second generator which contributes $t_1^2t_2$,
\item
the element (12)(34) which contributes $t_2^2$.
\end{itemize}
The cycle index is $(t_1^4+2t_1^2t_2+t_2^2)/4$.
This is replaced by the direct product
$(t_1^2+t_2)/2\times ({t'_1}^2+t'_2)/2$ as we wish to
represent the combined regions $\textit{1}\cup \textit{12}$ and 
$\textit{3}\cup \textit{23}$
that are related by one of the $C_2$ symmetries
differently from the regions $\overline{\textit 2}$
and $\underline{\textit 2}$ by the other $C_2$ symmetry.
\item
The previous diagram with a shrunk center circle, $^{3,5}\mathbb{X}_3$:
\begin{picture}(5,1.3)
\put(1.4,0.5){\circle{3.0}}
\put(2.0,0.5){\circle{0.5}}
\put(2.6,0.5){\circle{3.0}}
\put(2.8,0.5){1}
\put(1.9,0.5){2}
\put(1.0,0.5){3}
\end{picture}
The 7 regions inside the circles may be labeled by the circles that cover them: \textit{1}, \textit{12},
\textit{123}, \textit{23}, \textit{3}, \textit{$\overline{\textit{13}}$}, \textit{$\underline{\textit{13}}$},
using overline and underline to register the upper and lower regions.
The symmetry is the same as in the preceding diagram $^{3,4}\mathbb{X}_3$.
\item
The asymmetric bundle $^{3,6}\mathbb{X}_3$:
\begin{picture}(3,1.5)
\put(1.0,0.5){\circle{3.0}}
\put(1.1,0.5){\circle{0.8}}
\put(2.0,0.5){\circle{3.0}}
\put(2.0,0.5){1}
\put(0.5,0.5){2}
\put(1.1,0.5){3}
\end{picture}
The 5 regions inside the circles may be labeled by the circles that cover them: \textit{1}, \textit{12},
\textit{123}, \textit{23}, \textit{2}. The cycle index is $t_1$.
\end{enumerate}
Let $^3\mathbb{X}_N\supseteq {}^2\mathbb{X}_N$ denote the 
arrangements of $N$ nested circles which admit the topologies
of simple circles, the one topology of two-circle intersections and the
six topologies of three-circle intersections in the subregions.
\begin{exa}
This is an element of $^3\mathbb{X}_{10}^{(2)}$ which is not in $^2\mathbb{X}_N$:
\begin{picture}(3,1.5)
\put(1.0,0.5){\circle{3.0}}
\put(1.1,0.5){\circle{0.8}}
\put(1.0,0.5){\circle{0.4}}
\put(1.1,0.4){\circle{0.3}}
\put(2.0,0.5){\circle{3.0}}
\put(2.0,0.3){\circle{0.5}}
\put(2.2,0.6){\circle{0.5}}
\put(2.0,0.9){\circle{0.5}}
\put(3.0,0.8){\circle{0.6}}
\put(3.0,0.8){\circle{0.3}}
\end{picture}
\end{exa}

\subsection{Recurrences}
\begin{defn} Generating function of the topologies with up-to-three
intersections:
\begin{equation}
\sum_N |^3\mathbb{X}_N| z^N={}^3X(z).
\end{equation}
\end{defn}
The multiset interpretation as a forest of rooted trees with non-factoring
elements $^3\mathbb{X}_N^{(1)}$ in the roots holds again:
\begin{equation}
|^3\mathbb{X}_N^{(f)}| = \sum_{\pi(N): N=\{N_1^{c_1};N_2^{c_2};\ldots N_f^{c_f}\}} 
\prod_{j=1}^f
\binom{|^3\mathbb{X}_{N_j}^{(1)}|+c_j-1}{c_j}.
\end{equation}
There is one type of compound objects constructed by wrapping a circle
around others, one type of covering them with two intersecting circles,
and six types of covering them with three intersecting circles.
Because the two types $^{3,1}\mathbb{X}$
and
$^{3,2}\mathbb{X}$
of the 3-circles have the same number
of regions and the same symmetry, we count the first type twice
and drop the second;
because types $^{3,4}\mathbb{X}$ and and $^{3,5}\mathbb{X}$ 
have the same number of regions and the same
symmetry, we also count $^{3,4}\mathbb{X}$ twice and drop $^{3,5}\mathbb{X}$.
The upgrade of (\ref{eq.X2of1}) is
\begin{equation}
|^3\mathbb{X}_N^{(1)}| 
= |^3\mathbb{X}_{N-1}| +{}^2D_{N-2}
+2\,{}^{3,1}D_{N-3}
+{}^{3,3}D_{N-3}
+2\, {}^{3,4}D_{N-3}
+{}^{3,6}D_{N-3}.
\end{equation}
The generating functions are defined in the obvious way 
preserving the upper left type indices:
$\sum_{N\ge 0}{}\, ^{...}D_N z^N=\,{}^{...}D(z)$.
They are all anchored at $^{...}D_0=1$ and zero for negative $N$.
\begin{enumerate}
\item
The three regions in $^2\mathbb{X}_N$ are populated as before,
but now also accepting elements of $^3\mathbb{X}$ in their subregions
such that their values 
differ from the values of $\bar D_N$ of Equation (\ref{eq.Dbardef}):
\begin{equation}
^2D(z) = \frac12\, {}^3X(z)[^3X^2(z)+\,^3X(z^2)].
\end{equation}
\item
In $^{3,1}\mathbb{X}_N$ region \textit{123} is populated without
restriction. 
The remaining 6 regions associated via symmetry are then
incorporated with $t_j\mapsto{}\, ^3X^2(z^j)$, $j\ge 1$, in the cycle index, so
\begin{equation}
^{3,1}D(z) = \frac16\, ^3X(z)[^3X^6(z)+3\,{} ^3X^2(z)\, ^3X^2(z^2)+2\,{} ^3X^2(z^3)].
\end{equation}
\item
Region \textit{2} in $^{3,3}\mathbb{X}_N$ is populated without restrictions,
which contributes a factor $^3X(z)$. The pair regions $1$ and $12$
are populated without restrictions with is represented by $f(z)={}\,^3X^2(z)$.
The regions $3$ and $23$ associated to them via symmetry are then
incorporated with $t_j\mapsto f(z^j)$ in the cycle index, so
\begin{equation}
^{3,3}D(z) = \frac12\, ^3X(z)[^3X^4(z)+\,{}^3X^2(z^2)].
\end{equation}
\item
In $^{3,4}\mathbb{X}_N$ region \textit{123} is populated without
restriction which contributes a factor $^3X(z)$. 
Regions \textit{1} and \textit{12} are represented by $f(z)$
and $\overline{\textit{2}}$ is represented by $^3X(z)$.
Substituting $t_j=f(z^j)$, $t'_j ={}\,^3X(z)$ in the cycle index yields
\begin{equation}
^{3,4}D(z) = 
\frac14 \,^3X(z)[^3X^4(z)+\,{} ^3X^2(z^2)] [^3X^2(z)+\,{} ^3X(z^2)]
.
\end{equation}
\item
The five regions in $^{3,6}\mathbb{X}_N$ are populated without restrictions:
\begin{equation}
^{3,6}D(z) =\, ^3X^5(z).
\end{equation}
\end{enumerate}

\begin{table}
\begin{tabular}{r|r|rrrrrrrrrrrrrr}
$N$ & $ |^3\mathbb{X}_N|$ & 1 & 2 & 3 & 4 & 5 & 6 & 7 & 8 & 9 &\\
\hline
1 & 1& 1 & \\
2 & 3& 2 & 1 & \\
3 & 14& 11 & 2 & 1 & \\
4 & 61& 44 & 14 & 2 & 1 & \\
5 & 252& 169 & 66 & 14 & 2 & 1 & \\
6 & 1019& 609 & 323 & 70 & 14 & 2 & 1 & \\
7 & 4127& 2253 & 1431 & 356 & 70 & 14 & 2 & 1 & \\
8 & 17242& 8779 & 6320 & 1695 & 361 & 70 & 14 & 2 & 1 & \\
9 & 74007& 36319 & 27420 & 8081 & 1739 & 361 & 70 & 14 & 2 & 1 & \\
10 & 325615& 157297 & 119821 & 37849 & 8455 & 1745 & 361 & 70 & 14 & 2 & 1 & \\
11 & 1458604& 701901 & 528557 & 176894 & 40549 & 8510 & 1745 & 361 & 70 & 14 & 2 & 1 & \\
\end{tabular}
\caption{The number of topologies of nested $N$ circles intersecting at most
as triples, $|^3\mathbb{X}_N|$,
and the subcounts with $f$ factors, $|^3\mathbb{X}_N^{(f)}|$, $1\le f\le N$.
The row sums are the Euler transform of the column $f=1$.
}
\label{tab.X3}
\end{table}

The numerical evaluation of the recurrences leads to Table \ref{tab.X3}.
The first difference in comparison to Table \ref{tab.X2} is where
the 6 additional topologies offer new branches as soon
as at least 3 circles are involved: $|^3\mathbb{X}_3^{(1)}|
=|^2\mathbb{X}_3^{(1)}|+6$.

In a  wider context one would like to construct and count
all circle sets of $N$ circles with an arbitrary number of intersections
\cite[A250001]{EIS}. That is out of reach of this paper;
in Table \ref{tab.X3} that analysis is only complete up to $N=3$.

\bibliographystyle{amsplain}
\bibliography{all}

\end{document}